\documentclass[a4paper, 11pt, twoside, leqno]{article}
\usepackage{lipsum,mwe}
\usepackage[T1]{fontenc}
\usepackage[english]{babel}
\usepackage{lmodern}               
\usepackage{amsmath,amsthm,amssymb}
\usepackage{mathrsfs,esint,bbm,stmaryrd}
\usepackage{mathtools}
\usepackage{enumitem,graphicx,xcolor,subcaption}
\usepackage[textwidth=16cm,centering,headheight=24pt]{geometry}

\newtheorem{theorem}{Theorem}           

\newtheorem{proposition}[theorem]{Proposition}

\theoremstyle{definition}              
\newtheorem{definition}{Definition}

\theoremstyle{remark}                  

\DeclareMathOperator{\tr}{tr}                                       
\DeclareMathOperator{\curl}{curl}                                  
\let\div\relax
\DeclareMathOperator{\div}{div}                                     

\DeclareMathOperator{\SBV}{SBV}

\newcommand{\abs}[1]{\left| #1 \right|}                             

\DeclareMathAlphabet{\mathpzc}{OT1}{pzc}{m}{it}

\newcommand{\D}{\mathrm{D}}       
\renewcommand{\d}{\mathrm{d}}

\newcommand{\R}{\mathbb{R}}

\renewcommand{\SS}{\mathbb{S}}
\renewcommand{\P}{\mathbf{P}} 

\newcommand{\Q}{\mathbf{Q}}
\newcommand{\I}{\mathbf{I}}
\newcommand{\A}{\mathbf{A}}

\newcommand{\n}{\mathbf{n}}
\newcommand{\m}{\mathbf{m}}

\newcommand{\bfx}{\mathbf{x}}

\newcommand{\M}{\mathbf{M}}

\renewcommand{\u}{\mathbf{u}}
\newcommand{\NN}{\mathscr{N}}     

\renewcommand{\H}{\mathscr{H}}

\newcommand{\eps}{\varepsilon}
\newcommand{\bd}{{\mathrm{bd}}}

\SetSymbolFont{stmry}{bold}{U}{stmry}{m}{n}  
\newcommand{\nnu}{{\boldsymbol{\nu}}}

\renewcommand{\S}{\mathrm{S}}
\newcommand{\half}{\frac{1}{2}}
\newcommand{\Div}{{\rm div}\,}

\definecolor{lightblue}{rgb}{0.22,0.45,0.70}   
\definecolor{darkgray}{gray}{0.4}    
\definecolor{lightgray}{gray}{0.8}

\title{A free discontinuity model for smectic thin films}

\author{}

\author{John M.Ball\footnote{J. M. Ball,
       School of Mathematical \& Computer Sciences,
       Heriot-Watt University,
       Edinburgh, EH14 4AS, United Kingdom.
       E-mail: john.ball@hw.ac.uk}, \
 Giacomo Canevari\footnote{G. Canevari, 
        Dipartimento di Informatica, 
        Universit\`{a} di Verona, 
       Strada le Grazie 15, 37134 Verona, Italy.
        E-mail: giacomo.canevari@univr.it}, \
 Bianca Stroffolini\footnote{B. Stroffolini,
        Dipartimento di Matematica e Applicazioni,
        Universit\`{a} di Napoli Federico II, 
       Via Cintia, 80126 Napoli, Italy.
        E-mail: bstroffo@unina.it}
}

\date{\today}

\begin{document}

\maketitle

\begin{abstract}
We attempt to describe surface defects in smectic~A thin films by 
formulating a free discontinuity problem
--- that is, a variational problem in which the 
order parameter is allowed to have jump discontinuities
on some (unknown) set. The free energy functional contains 
an interfacial energy
which penalizes dislocations of the smectic layers
at the jump. We discuss mathematical issues related to
the existence of minimizers and provide examples
of minimizers in some simplified settings.

\end{abstract}

\section{Introduction}

\begin{figure}[h]
    \centering
    \begin{minipage}{0.25\textwidth}
        \centering
        \includegraphics[width=0.9\textwidth]{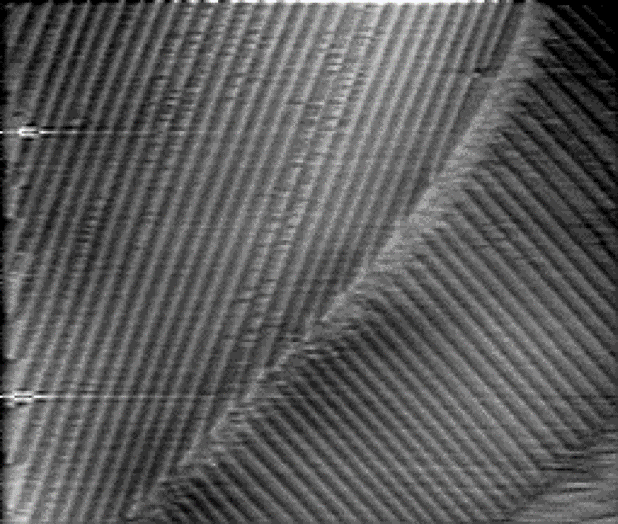} 
        \caption{AFM image of families of oily streaks of 8CB on ${\rm MoS}_2$ substrate.  Reprinted with permission from~\cite{michel2004} (copyright 2014, American Physical Society).}\label{smecticfig}
    \end{minipage}\hfill
    \begin{minipage}{0.7\textwidth}\vspace{-.37in}
        \centering
        \includegraphics[width=0.9\textwidth]{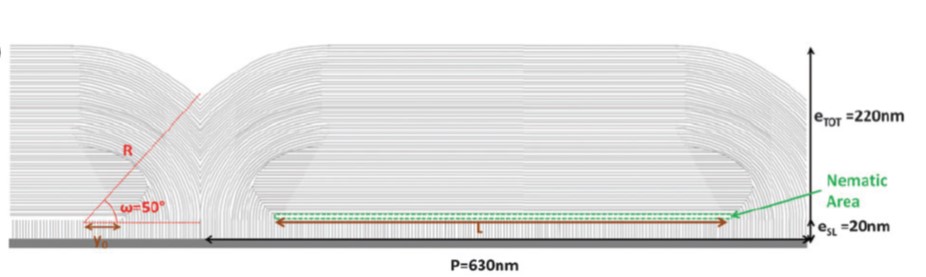} 
        \caption{Experimental representation of cross-sectional
layer configuration in an oily streak of 8CB on PVA substrate. Reproduced from~\cite{coursault2016} with permission from the Royal Society of Chemistry.}\label{smecticfig1}
    \end{minipage}
\end{figure}

In this paper we report on an attempt to formulate a mathematical model, capable of rigorous mathematical analysis, to describe aspects of the interesting experiments on smectic~A thin films of the group of  
E.~Lacaze~\cite{michel2004,Michel06, Zapponeetal2008, Zappone12,coursault2016}.
In these experiments thin films of 8CB liquid crystal are  deposited on various substrates, for example ${\rm MoS}_2$, mica and rubbed PVA. Depending on the thickness of the films and on the nature of the substrate different configurations of the smectic layers are observed~\cite{LacazeZappone}. The main features are illustrated in Figs~\ref{smecticfig}, \ref{smecticfig1}. Viewed from above (Fig.~\ref{smecticfig}) families of parallel `oily streaks' are observed, each of which consists of a flattened hemicylinder of smectic layers (Fig.~\ref{smecticfig1}). The representation in Fig.~\ref{smecticfig1} of the configuration of layers in an oily streak  on rubbed PVA substrate is not a direct observation, but is deduced from X-ray diffraction and ellipsometry; indeed recently some slight modifications to the likely configuration have been suggested by Niyonzima et al~\cite{niyonzima}.
The observed layer configurations are a consequence of the antagonistic boundary conditions, homeotropic on the upper surface in contact with air, and unidirectional planar anchoring on the bottom surface parallel to the substrate. Because for smectic A the director $\n$ is perpendicular to the layers, this means that the layers prefer energetically to be tangent  to the upper free surface and perpendicular to the substrate.

A feature of the observed layer configurations is the presence of surface (wall) defects across which the layer normal, and thus $\n$, jumps. This suggests that a good
mathematical framework to handle this problem is the free discontinuity setting originated by De Giorgi \& Ambrosio~\cite{degiorgiambrosio} for image segmentation and subsequently used by Francfort \& Marigo for Griffith fracture models~\cite{FrancfortMarigo}. The relevance of such a setting  for liquid crystal problems was proposed in~\cite{BallBedford}. The free discontinuity formulation uses a free energy in which there is a competition between bulk (volumetric) energy and interfacial energy corresponding to unknown surfaces of discontinuity. In the case of fracture mechanics the bulk energy is elastic energy and the interfacial energy is located on the unknown crack surfaces. For smectic A thin films we will take the bulk energy to be the elastic (Oseen-Frank) energy and the discontinuity surfaces will correspond to walls across which $\n$ jumps.

The thin film experiments described above would best be treated as a free boundary problem for a given volume of liquid crystal deposited on the substrate, with the aim of predicting the possible 3D configurations of oily streaks as well as their internal structure. A key prediction of such a model would be the width of the oily streaks. However we are not currently able to give conditions under which the minimum of the corresponding energy is attained and if so by what configurations of layers.

Instead we will work on a fixed domain $\Omega \subset \mathbb{R}^2$, to be thought of as the cross-section of an oily streak, and seek to predict the corresponding configurations of layers given suitable boundary conditions.

\section{The model}

\subsection{The elastic energy}
We use as a basic variable the director $\n$, $|\n|=1$,  which for smectic A is parallel to the layer normal $\m$, so that $\n\otimes\n=\m\otimes\m$. The three-dimensional Oseen-Frank energy density
\begin{equation} \label{b0}
 \begin{split}
  W(\n,\nabla\n)=\frac{1}{2}\left(K_1(\Div \n)^2\right.
   &+ K_2(\n\cdot\curl\n)^2 + K_3\abs{\n\wedge\curl\n}^2 \\
   &+(K_2+K_4)\left.(\tr(\nabla\n)^2-(\Div\n)^2)\right),
 \end{split}
\end{equation}
with Frank constants $K_i$, is invariant to changing the sign of $\n$, and thus can be expressed in terms of $\M:=\n\otimes\n$ and $\nabla\M$. For example, we have that 
\begin{equation} \label{b1}
 \begin{split}
  |\n\wedge\curl\n|^2 &= \abs{(\n\cdot\nabla)\n}^2 \\
  &=\half M_{rs}M_{ij,r}M_{ij,s},
 \end{split}
\end{equation}
where $M_{ij}=n_in_j$ and we have used $2n_in_{i,j}=(|\n|^2)_{,j}=0$.
We make the assumption of constant layer thickness (or locally parallel layers), which for sufficiently smooth $\m$ is equivalent to the condition $\curl \m=0$, and thus (assuming $\n$ sufficiently smooth so that we can take $\n=\m$) to the condition $\curl\n=0$.

We assume two-dimensional symmetry for~$\n$, so that 
 for~$\bfx=(x_1, x_2)\in \Omega$ 
\begin{equation}\n(\bfx) = ({\bf n}(x_1, x_2), 0)\label{b2}
\end{equation}
with ${\bf n}(x_1, x_2)\in\SS^1$, where~$\SS^1$ is the unit circle.
For such two-dimensional director fields,
an explicit computation shows that
the term $\tr (\nabla{\bf n} )^2 -(\div{\bf n})^2$ is identically equal to zero so that, taking into account also the constraint $\curl\n=0$, the Oseen-Frank energy \eqref{b0} reduces to
\begin{equation}
W({\bf n}, \nabla{\bf n})=\half K_1 |\nabla {\bf n}|^2,\label{b3}
\end{equation}
where we assume that $K_1>0$. This is the original model for smectics proposed by Oseen~\cite{Oseen33}. It can be viewed as a special case $\m=\n$ of the model of \cite{lesliestewartnakagawa}.
There are various other models for smectics
allowing variable layer thickness (\cite{deGennes72, chenlubensky, klemanparodi,hanetal}), that typically introduce the molecular number density~$\rho = \rho(\bfx)$ or fluctuations about it as a new macroscopic variable, with the smectic layers  being seen as density waves. These models can describe dislocations in smectic layers, for example,  although it is unclear how to understand the macroscopic variable $\rho$ varying over a molecular length-scale.

To allow for nonorientable configurations,
we reformulate the problem in terms of the two-dimensional $\Q$-tensor 
\begin{equation} \label{Q}
\Q:=\frac 1{\sqrt{2}} \left( {\bf n} \otimes {\bf n}  - \frac{\bf I} 2\right) 
\end{equation}
where~$\I$ is the~$2\times 2$ identity matrix.
The constant~$1/\sqrt{2}$ is a normalization factor,
which plays no essential role in our analysis.
At each point~$\bfx\in\Omega$, $\Q(\bfx)$ is a symmetric trace-free $2\times 2$ matrix which belongs to the set
\begin{equation} \label{N}
\NN:= \left \{ \frac 1{\sqrt{2}} \left( {\bf n} \otimes {\bf n}  - \frac{\bf I} 2\right)\colon {\bf n} \in \SS^1 \right \}.
\end{equation}
If~$\n\colon\Omega\to\SS^1$ is a smooth vector field and~$\Q$ is defined as in~\eqref{Q}, an explicit computation shows that~$\abs{\nabla{\bf n}}^2 
= \abs{\nabla \Q}^2 := Q_{ij,k} \, Q_{ij,k}$,
expressing the elastic energy \eqref{b3} in terms of~$\Q$. 

We also need to express the constraint~$\curl\n = 0$ in terms of~$\Q$.
We define
\begin{equation*}
{\bf A}(\Q)(\nabla \Q, \nabla \Q) 
:=\left( \sqrt{2} Q_{hk}+\frac12 \delta_{hk}\right) Q_{ij,h} \, Q_{ij,k}
\end{equation*}
This is a quadratic form in~$\nabla\Q$, reminiscent of 
the cubic term in the Landau-de Gennes elastic energy (see e.g. \cite[Section 4]{BallILCC}). 
Thanks to~\eqref{b1}, \eqref{Q}, we find that 
\begin{equation*}
 \begin{split}
  {\bf A}(\Q)( \nabla \Q, \nabla \Q) = |\n\wedge\curl\n|^2,
 \end{split}
\end{equation*}
and since $\abs{\curl\n}^2=\abs{\n\wedge\curl\n}^2+(\n\cdot\curl\n)^2$, and  keeping in mind that a two-dimensional vector field
is orthogonal to its curl, we obtain 
$\A(\Q)(\nabla\Q, \, \nabla\Q) = \abs{\curl\n}^2$. 
Therefore,  we impose the constraint
\begin{equation} \label{A0}
 \A(\Q)(\nabla\Q, \, \nabla\Q) = 0,
\end{equation}
which expresses constant layer thickness in terms of $\Q$. 

\subsection{The function space}

In order to specify in precise mathematical terms any energy minimization problem, it is necessary to say in which function space the minimum is sought. The function space describes the allowed singularities of the unknown function or map and is part of the model. Making the function space larger, so that worse singularities are allowed, may lead to different  minimizers (this is known as the  {\it Lavrentiev phenomenon}, see \cite[Section 3]{BallILCC}). The main function space used for free discontinuity problems, and developed by De Giorgi and Ambrosio~\cite{degiorgiambrosio}, is the space SBV of {\it special functions of bounded variation}.
In fact, in this work we will consider a slight variant of the space~$\SBV$,
i.e. the space~$\SBV^2$. Let~$\Omega$ be an open, bounded planar region.
In technical terms, a map~$\u\colon\Omega\to\R^m$
belongs to~$\SBV^2(\Omega, \, \R^m)$ if its distributional derivative~$\D\u$
is a finite measure with no Cantor part, and whose
absolutely continuous part~$\nabla\u$ is square integrable
(see e.g~\cite{degiorgiambrosio, AmbrosioFuscoPallara}
for more details). For the purposes of this paper,
the key points of the definition and theory are the following:
\begin{enumerate}
\item{} for any ${\bf u}\in \SBV^2(\Omega, \, \R^m)$ 
there is a one-dimensional {\it jump
set} $\S_{\bf u}$ consisting of the points $\bfx\in \Omega$ at which ${\bf u}$ has a jump discontinuity,
\item{} for (almost) any point $\bfx\in\S_{\bf u}$ there is a well-defined unit normal $\nnu(\bfx) = \nnu_{\u}(\bfx)$ to $\S_{\bf u}$ and well-defined limits ${\bf u}^+(\bfx), {\bf u}^-(\bfx)$ from either side of $\S_{\bf u}$. There may be an exceptional
set of points~$\bfx\in\S_{\u}$ at which~$\nnu(\bfx)$, $\u^+(\bfx)$ or~$\u^-(\bfx)$
are not defined, but this must be a set of zero length;
\item{} the gradient $\nabla {\bf u}$ is defined in $\Omega \setminus\S_{\bf u}$
and $\displaystyle\int_{\Omega}|\nabla {\bf u}(\bfx)|^2 \, \d\bfx < +\infty$.
\end{enumerate}
The space~$\SBV^2(\Omega, \, \R^{2\times 2})$, whose elements
are~$2\times 2$~matrix fields on~$\Omega$ in the class~$\SBV^2$,
is defined analogously by identifying $\R^{2\times 2}$ with $\R^4$.

%

As we are considering a free-discontinuity problem,
the map~$\Q$ is allowed to have  jump discontinuities 
on a one-dimensional set~$\S_{\Q}$, and the condition~\eqref{A0}
loses its meaning at points of~$\S_{\Q}$. Therefore, the
constraint~\eqref{A0} is only enforced in the complement~$\Omega\setminus\S_{\Q}$.
In technical terms,  \eqref{A0} only involves
the absolutely continuous part~$\nabla\Q$ of the
distributional derivative of~$\Q$.

\subsection{The jump energy}

\begin{figure}[tb]
	\centering
	 \includegraphics[height=.19\textheight]{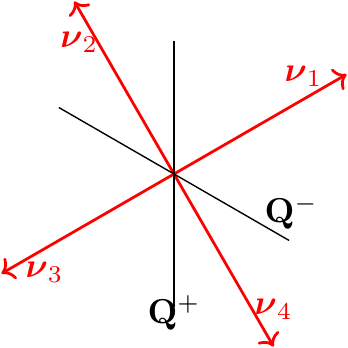}
	 \hspace{.3cm}
     \includegraphics[height=.19\textheight]{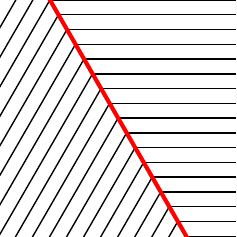}
     \hspace{.3cm}
     \includegraphics[height=.19\textheight]{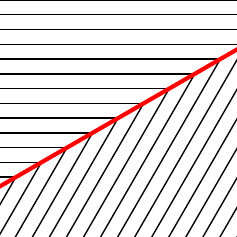}
	\caption{Left: given~$\Q^+$ and~$\Q^-$, there are four 
	unit vectors~$\nnu_1$, $\nnu_2$, $\nnu_3$, $\nnu_4$
	that satisfy~\eqref{nodisl-n}.
	The~$\Q$-tensors~$\Q^+$ and~$\Q^-$
	are represented by (arrowless black) straight lines,
	oriented parallel to the molecular directors~$\n^+$, 
	$\n^-$ corresponding to~$\Q^+$, $\Q^-$ respectively.
	Here~$\Q^+$ and $\n^+$ (respectively, $\Q^-$ and~$\n^-$)
	are related  to each other by~\eqref{Q}.
	Centre and right: two piecewise constant 
	configurations that satisfy~\eqref{nodisl-n}. The black lines 
	represent the smectic layers, which are orthogonal 
	to the molecular director, and the thick red line represents the jump.}
	\label{fig:bisectors}
\end{figure}

Let
\begin{equation} \label{SBVN}
 \SBV^2(\Omega, \, \NN) :=\left\{\Q\in\SBV^2(\Omega, \, \R^{2\times 2}),
 \ \Q(x)\in\NN \ \textrm{ for } x\in\Omega\right\}
\end{equation}
be the set of~$\NN$-valued $\Q$-tensors in the class~$\SBV^2$.
For~$\Q\in\SBV^2(\Omega, \, \NN)$, we consider the jump
energy
\begin{equation} \label{jumpenergy}
 \int_{\S_{\Q}}\varphi(\Q^+, \, \Q^-, \, \nnu) \, \d\H^1
\end{equation}
where the integral is taken over the (one-dimensional) jump set~$\S_{\Q}$,
with respect to the length measure~$\d\H^1$, 
and~$\varphi\colon\NN\times\NN\times\SS^1\to\R$ is a continuous function.
We want to choose the jump energy density~$\varphi$ so that~\eqref{jumpenergy}
is a good model for the energy of a defect wall. 
Natural conditions on~$\varphi$ are:
\begin{enumerate}[label=(C\arabic*), ref=C\arabic*]
 \item \label{C:first} $\varphi(\Q^+, \, \Q^-, \, \nnu) = 0$ if~$\Q^+ = \Q^-$;
 \item invariance with respect to the orientation of the jump set, that is
 \[
  \varphi(\Q^-, \, \Q^+, \, -\nnu) = \varphi(\Q^+, \, \Q^-, \, \nnu)
 \]
 for any~$\Q^+\in\NN$, $\Q^-\in\NN$, $\nnu\in\SS^1$;
 \item \label{C:rot} frame-indifference, that is
 \[
  \varphi(\mathbf{R}\Q^+\mathbf{R}^{\mathsf{T}}, \, \mathbf{R}\Q^-\mathbf{R}^{\mathsf{T}}, \mathbf{R}\nnu)
  = \varphi(\Q^+, \, \Q^-, \, \nnu)
 \]
 for all~$\Q^+\in\NN$, $\Q^-\in\NN$, $\nnu\in\SS^1$ and 
 all  orthogonal matrices~$\mathbf{R}\in\mathrm{O}(2)$.
\end{enumerate}
Another condition we impose is that $\varphi$ should penalize
dislocations of the smectic layers. Generically, it may
not be possible to define the smectic layers consistently across the jump set,
because there may be dislocations. A geometric condition 
for the layers to match at the jump set
is that, at each jump point, the normal to the jump set~$\nnu$
bisects any of the angles between the smectic layers.
For smectic~A liquid crystals, the layers are orthogonal to
the  the molecular directors~$\n^+$, $\n^-$ on either side of the jump,
so bisecting any of the angles between the layers is equivalent
to bisecting any of the angles between~$\n^+$ and~$\n^-$.
This can be written as
\begin{equation} \label{nodisl-n}
 (\n^+\cdot\nnu)^2 = (\n^-\cdot\nnu)^2
\end{equation}
or equivalently, in terms of the~$\Q$-tensor, as
\begin{equation} \label{nodisl}
 \Q^+\nnu\cdot\nnu = \Q^-\nnu\cdot\nnu.
\end{equation}
In order to penalize dislocations,
we therefore impose the following condition on~$\varphi$:
\begin{enumerate}
[label=(C\arabic*), ref=C\arabic*, resume]
 \item \label{C:disl} for given~$\Q^+, \Q^-\in\NN$, the function~$\nnu\in\SS^1\mapsto\varphi(\Q^+, \, \Q^-, \, \nnu)$
 is minimized for a~$\nnu$ satisfying~\eqref{nodisl}.
\end{enumerate}
Given distinct~$\Q^+$ and~$\Q^-$, there are four such~$\nnu$ that bisect
the angles between the corresponding directors~$\n^+$, $\n^-$;
we label them as~$\nnu_1$, $\nnu_2$, $\nnu_3 = -\nnu_1$ and~$\nnu_4 = -\nnu_2$,
as in Fig.~\ref{fig:bisectors}. 
By taking~$\mathbf{R} = -\I$ in condition~\eqref{C:rot},
we obtain that~$\varphi(\Q^+, \, \Q^-, \, \nnu_1) = \varphi(\Q^+, \, \Q^-, \, \nnu_3)$ and~$\varphi(\Q^+, \, \Q^-, \, \nnu_2) = \varphi(\Q^+, \, \Q^-, \, \nnu_4)$. However, we must account for the possibility that, in general,
$\varphi(\Q^+, \, \Q^-, \, \nnu_1)\neq\varphi(\Q^+, \, \Q^-, \, \nnu_2)$. 
Indeed, configurations whose jump sets are oriented according to the unit normal
$\nnu_1$ or~$\nnu_2$ have (generically) different geometric properties,
because~$\Q^+\nnu_1\cdot\nnu_1 \neq \Q^+\nnu_2\cdot\nnu_2$
unless the directors~$\n^+$, $\n^-$ corresponding to~$\Q^+$, $\Q^-$
are orthogonal to each other. The condition~\eqref{C:disl}
is compatible with an energy density that, as a function of~$\nnu$,
has two global minima at~$\nnu_1$, $\nnu_3$ (or, respectively,
at~$\nnu_2$, $\nnu_4$) and, say, two local minima 
at~$\nnu_2$, $\nnu_4$ (respectively, at~$\nnu_1$, $\nnu_3$),
at a possibly higher energy value.

\begin{figure}[tb]
	\centering
     \includegraphics[width=.44\textwidth]{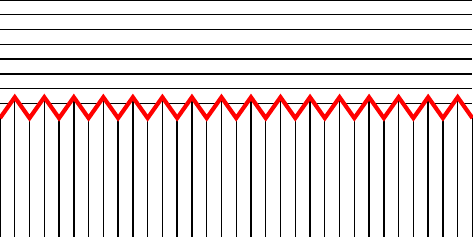}
     \hspace{.5cm}
     \includegraphics[width=.44\textwidth]{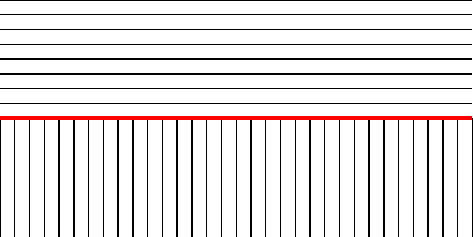}
	\caption{A mathematical pathology associated with
	the energy density~$\zeta_\alpha$, given by~\eqref{zeta_alpha}.
	The black lines represent the smectic layers, while the
	thick red line represents the jump set.
	On the left, a piecewise constant configuration in a rectangle,
	with a zig-zag interface
	and zig-zag angles of~$45$ degrees; on the right, 
	a piecewise constant configuration with a horizontal 
	jump set. If we choose the energy density as in~\eqref{zeta_alpha},
	the energy of the configuration on the left is~$\sqrt{2}b\mu$,
	where~$b$ is the length of the bottom base of the rectangle,
	irrespective of the spacing of the zig-zag.
	As the latter tends to zero,
	the configuration on the left converges (in a suitable sense)
	to that on the right, which costs an infinite energy.}
	\label{fig:zigzag}
\end{figure}

A singular jump energy satisfying all the conditions~\eqref{C:first}--\eqref{C:disl} is given by
\begin{equation} \label{zeta_alpha}
 \zeta_\alpha(\Q^+, \, \Q^-, \, \nnu) :=
 \begin{cases}
  \abs{\Q^+ - \Q^-}^\alpha 
  &\textrm{if } \Q^+\nnu\cdot\nnu = \Q^-\nnu\cdot\nnu \\
  +\infty &\textrm{otherwise,}
 \end{cases}
\end{equation}
where~$\alpha$ is a parameter such that~$0 < \alpha < 1$.
Choosing~$\alpha > 0$ guarantees that the condition~\eqref{C:first}
is satisfied. On the other hand, taking~$\alpha < 1$ is required 
for reasons of mathematical consistency;
had we taken~$\alpha \geq 1$, then there would be
no guarantee that a minimizer for the energy functional exists 
in the space~$\SBV^2(\Omega, \, \NN)$. 
(We do not discuss this issue in detail
and refer to e.g.~\cite[Sections~4.1--2]{AmbrosioFuscoPallara}.)
However, even if the parameter~$\alpha$ is chosen carefully,
the functional associated with~\eqref{zeta_alpha}
suffers from a mathematical pathology, which is illustrated
in Fig.~\ref{fig:zigzag}: there exist sequences~$\Q_j\in\SBV(\Omega, \, \NN)$
that converge to a limit map~$\Q\in\SBV(\Omega, \, \NN)$ in a suitable sense,
yet the energy of~$\Q$ is strictly larger than the 
limit, as~$j\to+\infty$, of the energy of~$\Q_j$.
This behaviour, known as lack of lower semicontinuity,
is particularly evident in Fig.~\ref{fig:zigzag}, because
the energy of the limit configuration is infinite.
However, even if we considered a modified energy density 
that takes only finite values, this pathological behaviour 
could persist. The main issue is that the energy density~$\zeta_\alpha$
is not \emph{BV-elliptic}.
BV-ellipticity, which was introduced by Ambrosio and Braides~\cite{AmbrosioBraides-I},
is a necessary condition for the
lower semicontinuity of the energy functional,
and is an important assumption to ensure the existence of minimizers of
free-discontinuity problems~\cite{AmbrosioBraides-I, Ambrosio90}.

\begin{figure}[tb]
	\centering
	 \includegraphics[height=.3\textheight]{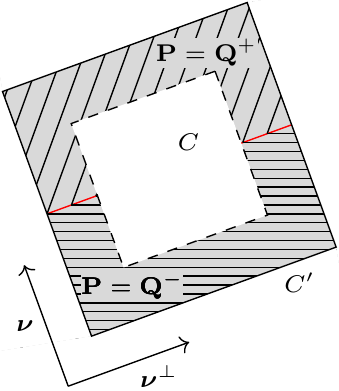}
	\caption{The domain and boundary conditions
	in the definition of BV-ellipticity (Definition~\ref{def:BVelliptic}).
	The domain~$C$ is a unit square, rotated in such a way 
	that the sides are parallel to~$\nnu$, $\nnu^\perp$.
	The `boundary conditions', defined in a collar of~$\partial C$,
	are piecewise constant and are defined by~$\Q^+$ and~$\Q^-$. 
	The admissible configurations~$\P$ are piecewise constant inside~$C$
	(and are allowed to take finitely many values only).}
	\label{fig:BVelliptic}
\end{figure}

We briefly explain the notion of BV-ellipticity.
Let~$\Q^+\in\NN$, $\Q^-\in\NN$, $\nnu\in\SS^1$ be given.
Let~$C$ be a (closed) unit square, centred at the origin,
whose sides are parallel to~$\nnu$ and to~$\nnu^\perp := (\nu_2, \, -\nu_1)$
(see Figure~\ref{fig:BVelliptic}).
Let~$C^\prime$ be a slightly larger such square,
such that~$C$ is contained in the interior of~$C^\prime$.
We consider a piecewise constant map~$\Q^*\in\SBV(C^\prime, \, \NN)$, defined as
\begin{equation} \label{Q*}
 \Q^*(\bfx) := \begin{cases} 
              \Q^+ & \textrm{if } \bfx\cdot\nnu > 0 \\
              \Q^- & \textrm{if } \bfx\cdot\nnu < 0
             \end{cases}
\end{equation}
The map~$\Q^*$ jumps along a straight line orthogonal to~$\nnu$.
The notion of BV-ellipticity is defined in terms of a suitable
minimization problem. Let~$\mathcal{C}[\Q^+, \, \Q^-, \, \nnu]$
be the class of all maps~$\P\in\SBV^2(C^\prime, \, \NN)$ that satisfy
the following properties: 
\begin{enumerate}[label=(\roman*)]
 \item $\P = \Q^*$ in~$C^\prime\setminus C$;
 \item $\P$ has finite range (i.e., $\P$ is a piecewise 
 constant configuration that only takes finitely many values).
\end{enumerate}
The condition~(i) plays the r\^ole of a boundary condition.
As for condition~(ii), restricting our attention to piecewise constant
configurations with finite range allows us to neglect
all elastic contributions for the moment and focus on the jump energy.

\begin{definition}[BV-ellipticity] \label{def:BVelliptic}
 We say that a function~$\varphi\colon\NN\times\NN\times\SS^1\to\R$
 is BV-{\it elliptic} if and only if, for any~$(\Q^+, \, \Q^-, \, \nnu)\in\NN\times\NN\times\SS^1$, we have
 \[
  \inf_{\P\in\mathcal{C}[\Q^+, \, \Q^-, \, \nnu]}
  \int_{\S_{\P}\cap C} \varphi(\P^+, \, \P^-, \, \nnu_{\P}) \, \d\H^1
  = \varphi(\Q^+, \, \Q^-, \, \nnu)
 \]
\end{definition}

Equivalently, $\varphi$ is BV-elliptic if and only if,
for any~$(\Q^+, \, \Q^-, \, \nnu)\in\NN\times\NN\times\SS^1$,
the map~$\Q^*$ defined in~\eqref{Q*}
minimizes the jump energy functional associated with~$\varphi$
among all competitors in~$\mathcal{C}[\Q^+, \, \Q^-, \, \nnu]$.
Roughly speaking, a jump energy density~$\varphi$
is BV-elliptic if the corresponding jump energy functional
favours the jump set to be (locally) a straight line.
This indeed rules out pathological behaviour,
such as the one discussed in Fig.~\ref{fig:zigzag}
(see e.g.~\cite{AmbrosioBraides-I, AmbrosioFuscoPallara}).
In practice, deciding whether a given function is BV-elliptic
or not may not an easy task.
Sufficient and necessary conditions for BV-ellipticity
have been proposed in the literature
(see e.g.~\cite{AmbrosioBraides-I,  AmbrosioBraidesII-90, AmbrosioFuscoPallara, Caraballo13}), but even these might not 
be immediately applicable to concrete examples.

The function~$\zeta_\alpha$ is not BV-elliptic,
because it fails to satisfy some necessary conditions
for BV-ellipticity (such as convexity in the variable~$\nnu$
--- see e.g.~\cite[Theorems~5.11 and~5.14]{AmbrosioFuscoPallara}). 
However, a natural candidate
for a BV-elliptic function that satisfies~\eqref{C:first}--\eqref{C:disl}
is the \emph{BV-elliptic envelope} of~$\zeta_\alpha$,
that is, the largest BV-elliptic function~$\varphi$ such 
that~$\varphi\leq \zeta_\alpha$.
As it turns out, the BV-elliptic envelope of~$\zeta_\alpha$
can be calculated explicitly, and is given by
\begin{equation} \label{phi_alpha}
 \varphi_\alpha(\Q^+, \, \Q^-, \, \nnu)
 := \abs{\Q^+ - \Q^-}^\alpha \left(1 + \frac{\sqrt{2}\abs{\Q^+\nnu\cdot\nnu - \Q^-\nnu\cdot\nnu}}{\abs{\Q^+ - \Q^-}}\right)^{1/2} 
\end{equation}
if~$\Q^+\neq\Q^-$, and~$\varphi_\alpha(\Q^+, \, \Q^-, \, \nnu) = 0$
if~$\Q^+ = \Q^-$. The function~$\varphi_\alpha$ can be
expressed in terms of angular variables, which is sometimes
convenient for computations. If~$\n^+ = (\cos\beta_+, \, \sin\beta_+)$,
$\n^- = (\cos\beta_-, \, \sin\beta_-)$ are the directors
corresponding to~$\Q^+$, $\Q^-$ respectively and
$\nnu= (\cos\gamma, \, \sin\gamma)$, with~$\beta_+$, $\beta_-$, $\gamma$
real numbers, then
\begin{equation} \label{phi_alpha_bis}
 \begin{split}
  \varphi_\alpha(\Q^+, \, \Q^-, \, \nnu)
  &= \abs{\sin(\beta_+ - \beta_-)}^\alpha \left(1 + \abs{\sin\left(\beta_+ + \beta_- - 2\gamma\right)}\right)^{1/2} \\
  &= \abs{\sin(\beta_+ - \beta_-)}^\alpha \left(\abs{\cos\left(\frac{\beta_+ + \beta_-}{2} - \gamma\right)} + \abs{\sin\left(\frac{\beta_+ + \beta_-}{2} - \gamma\right)}\right)
 \end{split}
\end{equation}
Equation~\eqref{phi_alpha_bis} is obtained from~\eqref{phi_alpha},
by substituting~\eqref{Q} for~$\Q^+$, $\Q^-$ and applying
trigonometric identities.

The proof that~$\varphi_\alpha$
is indeed BV-elliptic (and, in fact, is the BV-elliptic envelope of~$\zeta_\alpha$)
is rather technical, and { we will present it in a
forthcoming paper~\cite{BCS}\footnote{{ In~\cite{BCS},
we will show that~$\varphi_\alpha$
is not only BV-elliptic, but also jointly convex, in the sense of~\cite[Definition~5.17]{AmbrosioFuscoPallara}.
Joint convexity is a sufficient condition for BV-ellipticity, but it is 
not known whether it is a necessary condition as well.}}.}
However, from~\eqref{phi_alpha} and~\eqref{phi_alpha_bis}, it is immediate to
check that~$\varphi_\alpha$ satisfies~\eqref{C:first}--\eqref{C:disl}.
For given values~$\Q^+\neq\Q^-$, the function~$\nnu\mapsto\varphi_\alpha(\Q^+, \, \Q^-, \, \nnu)$ has four global minima, at the same energy value,
corresponding exactly to the four unit vectors~$\nnu$ that satisfy~\eqref{nodisl}
(i.e. $\gamma = (\beta_+ + \beta_-)/2 + k\pi/2$, for integer~$k$).
We do not have an example of an energy density that satisfies
\eqref{C:first}--\eqref{C:disl}, is BV-elliptic, and has non-equal
local minima at the four vectors~$\nnu$ that satisfy~\eqref{nodisl}.

\section{Existence of minimizers and examples}
\label{sect:results}

When the jump energy density is chosen as in~\eqref{phi_alpha},
we can prove existence of minimizers for the energy functional,
subject to appropriate boundary conditions. 
For instance, strong anchoring at the boundary is 
often represented mathematically by imposing
Dirichlet boundary conditions,
of the form~$\Q = \Q_{\bd}$ on~$\partial\Omega$,
where~$\Q_{\bd}\colon\partial\Omega\to\NN$ is a boundary datum. 
However, these Dirichlet boundary conditions are not suited
to free-discontinuity problems. Instead, one has to define the boundary conditions carefully in order to cater for the possibility that $\Q$ jumps at the boundary. This is done by \lq thickening the boundary\rq. In other words,
we consider an exterior {{neighbourhood}}~$\Gamma$ of the 
boundary~$\partial\Omega$ and impose the condition
\begin{equation} \label{bc}
 \Q = \Q_\bd \qquad \textrm{on } \Gamma,
\end{equation}
where the datum~$\Q_{\bd}$ is now defined in~$\Gamma$.
The map~$\Q_{\bd}$ needs to be compatible with our setting 
(i.e., $\Q_{\bd}\in\SBV^2(\Gamma, \, \NN)$ and it must satisfy~$\A(\Q_{\bd})(\nabla\Q_{\bd}, \, \nabla\Q_{\bd}) = 0$ in~$\Gamma$).

Using the direct methods in the Calculus of Variations, 
we could prove the following existence theorem:

\begin{theorem} \label{th:existence}
For any $K_1 > 0$, $\mu>0$ and~$0 < \alpha < 1$, the energy functional:
\begin{equation} \label{energy}
I(\Q)= \frac{K_1}2 \int_{\Omega} |\nabla \Q|^2 \d x +\mu \int_{\S_{\Q}} \varphi_{\alpha} (\Q^+, \Q^-, {\bf \nu}) \, \d \H^1
\end{equation}
with~$\varphi_\alpha$ as in~\eqref{phi_alpha}, 
attains a minimum among $\Q\in\SBV^2 (\Omega, \NN)$ satisfying ${\bf A}(\Q)( \nabla \Q, \nabla \Q)=0$  and the boundary conditions~\eqref{bc}.
\end{theorem}

The proof of Theorem~\ref{th:existence} will be given in~\cite{BCS}.
To try to get some understanding of the behaviour of minimizers,
we consider a simplified (or over-simplified) problem,
for which we can find the minimizer explicitly.
We consider a rectangular domain, $\Omega = (-L, \, L)\times (0, \, H)$,
with~$L > 0$ and~$H > L/2$. We focus  attention
on a restricted class of configurations, 
whose jump set can be described in polar coordinates as the curve
\[
 \mathcal{C} = \left\{(\rho(\theta)\cos\theta, \, \rho(\theta)\sin\theta)\colon 0 \leq \theta\leq \pi\right\}
\]
where~$\rho\colon [0, \, \pi]\to\R$ is a scalar function
to be determined, subject to the conditions~$\rho(0) = \rho(\pi) = L$.
Moreover, we assume that  $\Q$ is given by
\begin{equation} \label{Qrho}
 \Q(\bfx) := \frac{1}{\sqrt{2}}
 \left(\n(\bfx)\otimes\n(\bfx) - \frac{\I}{2}\right)
\end{equation}
where the director~$\n$ is defined in polar coordinates as
\[
 \n(r\cos\theta, \, r\sin\theta) 
 := \begin{cases}
  (\cos\theta, \, \sin\theta) & \textrm{if } r < \rho(\theta)\\
  (0, \, 1) &\textrm{otherwise.}
 \end{cases}
\]
This configuration~$\Q$ is uniquely determined by
the function~$\rho$, $\Q = \Q[\rho]$; it has
horizontal layers above the curve~$\mathcal{C}$ and concentric circular
layers below~$\mathcal{C}$. In particular, $\Q$ satisfies
planar anchoring conditions at the boundary~$x_2 = 0$,
homeotropic anchoring conditions at~$x_2 = H$
and periodic boundary conditions at~$x_1 = \pm L$.
To simplify the problem further, we neglect the elastic energy.
(Heuristically, we expect this approximation to
be relevant in the limit as~$\mu/K_1\to+\infty$.)

\begin{figure}[tb]
  \centering
  \includegraphics[width=.9\textwidth]{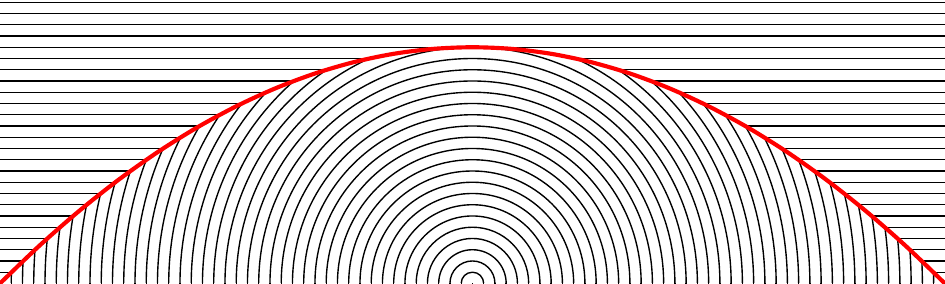}
  \caption{The minimizing configuration given by Proposition~\ref{prop:parabola}.
  The black lines represent the smectic layers, while
  the thick red line is the jump set~$\mathcal{C}$.}
  \label{fig:parabola}
\end{figure}

\begin{proposition} \label{prop:parabola}
 For any~$0 < \alpha < 1$, the unique minimizer of the functional
 \[
  \rho\mapsto \int_{\mathcal{C}} \varphi_\alpha(\Q^+[\rho], \, \Q^-[\rho], \, \nnu_{\Q[\rho]}) \, \d\H^1,
 \]
 among all $($Lipschitz continuous$)$ functions~$\rho\colon[0, \, \pi]\to\R$
 such that~$\rho(0) = \rho(\pi) = L$, is given by
 \begin{equation} \label{parabola-rho}
  \rho(\theta) := \frac{L}{1 + \sin\theta}
 \end{equation}
 for~$\theta\in [0, \, \pi]$.
\end{proposition}

The proof of Proposition~\ref{prop:parabola} will be given in~\cite{BCS}.
When~$\rho$ is given by~\eqref{parabola-rho},
the jump set~$\mathcal{C}$ can be described in Cartesian coordinates
as the graph of
\[
 x_2 = \frac{L}{2} - \frac{x_1^2}{2L}
 \qquad \textrm{for } -L \leq x_1 \leq L.
\]
In particular, $\mathcal{C}$ is a parabolic arc (see Fig.~\ref{fig:parabola}).
It can be shown that, at each point, the normal to the curve~$\mathcal{C}$
bisects the angle between the smectic layers.

Next we consider a domain which is a quarter disk of unit radius,
$\Omega := \{\bfx = (x_1, \, x_2)\in\R^2\colon x_2^2 + x_2^2 < 1, \ x_1 > 0, \ x_2 > 0\}$. Again, we minimize within a restricted class of configurations,
whose jump set has the form
\begin{equation} \label{C-num}
 \mathcal{C} = \left\{(\rho(\theta)\cos\theta, \, \rho(\theta)\sin\theta)\colon 0 \leq \theta\leq \frac{\pi}{2}\right\}
\end{equation}
The function~$\rho$ is unknown but, in contrast to 
the previous example, the boundary values~$\rho(0)$, 
$\rho(\pi/2)$ are unspecified.
We consider maps~$\Q$ of the form~\eqref{Qrho}, 
where the director~$\n$ is given by
\begin{equation} \label{director-num}
 \n(r\cos\theta, \, r\sin\theta) := 
 \begin{cases}
  (0, \, 1) & \textrm{if } r < \rho(\theta), \ 0 < \theta < \pi/2\\
  (\cos\theta, \, \sin\theta) & \textrm{if }
   r > \rho(\theta), \ 0 < \theta < \pi/2.
 \end{cases}
\end{equation}
In particular, the layers are horizontal on the left side of~$\mathcal{C}$
and circular on the right side. Moreover,
we impose planar anchoring conditions near
the boundary at~$x_1 = 0$, defined in terms 
of the director as~$\n(x_1, \, x_2) = (1, \, 0)$ for~$x_2 < 0$.
As a result, the jump set of~$\Q$ contains an additional
component on the boundary of~$\partial\Omega$,
i.e. the straight line segment of endpoints $(0, \, 0)$
and~$(0, \, \rho(0))$, where the director jumps from 
the value~$(1, \, 0)$ to~$(0, 1)$. In the energy functional,
we account for this additional contribution as well.
Moreover, in this example
we do \emph{not} neglect the elastic energy.

\begin{figure}[htb]
  \centering
  \begin{minipage}{.49\textwidth}
   \centering
   \includegraphics[width=.98\textwidth]{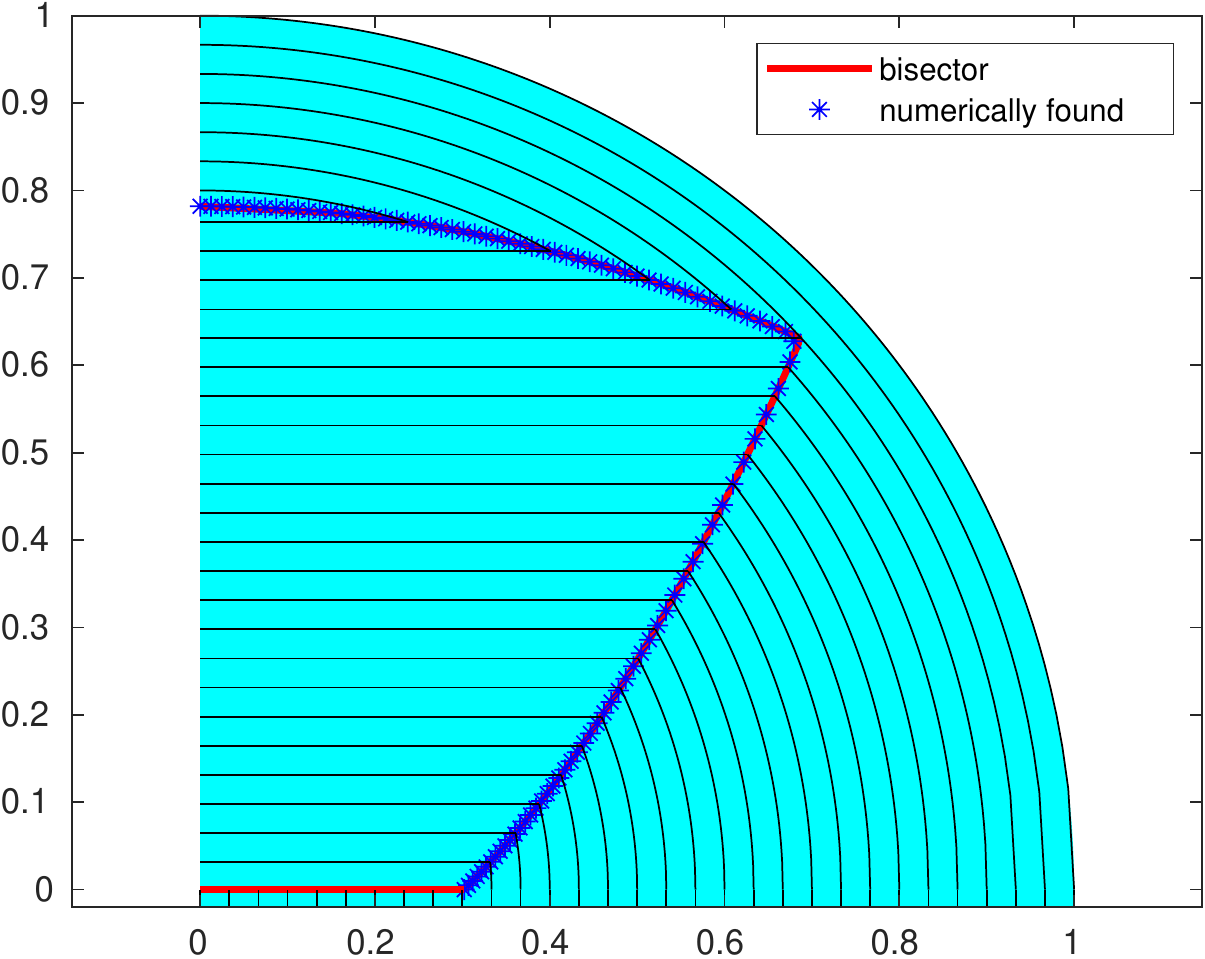}
   
   $\mu = 1$, $\alpha = 0.5$
  \end{minipage} \
  \begin{minipage}{.49\textwidth}
   \centering
   \includegraphics[width=.98\textwidth]{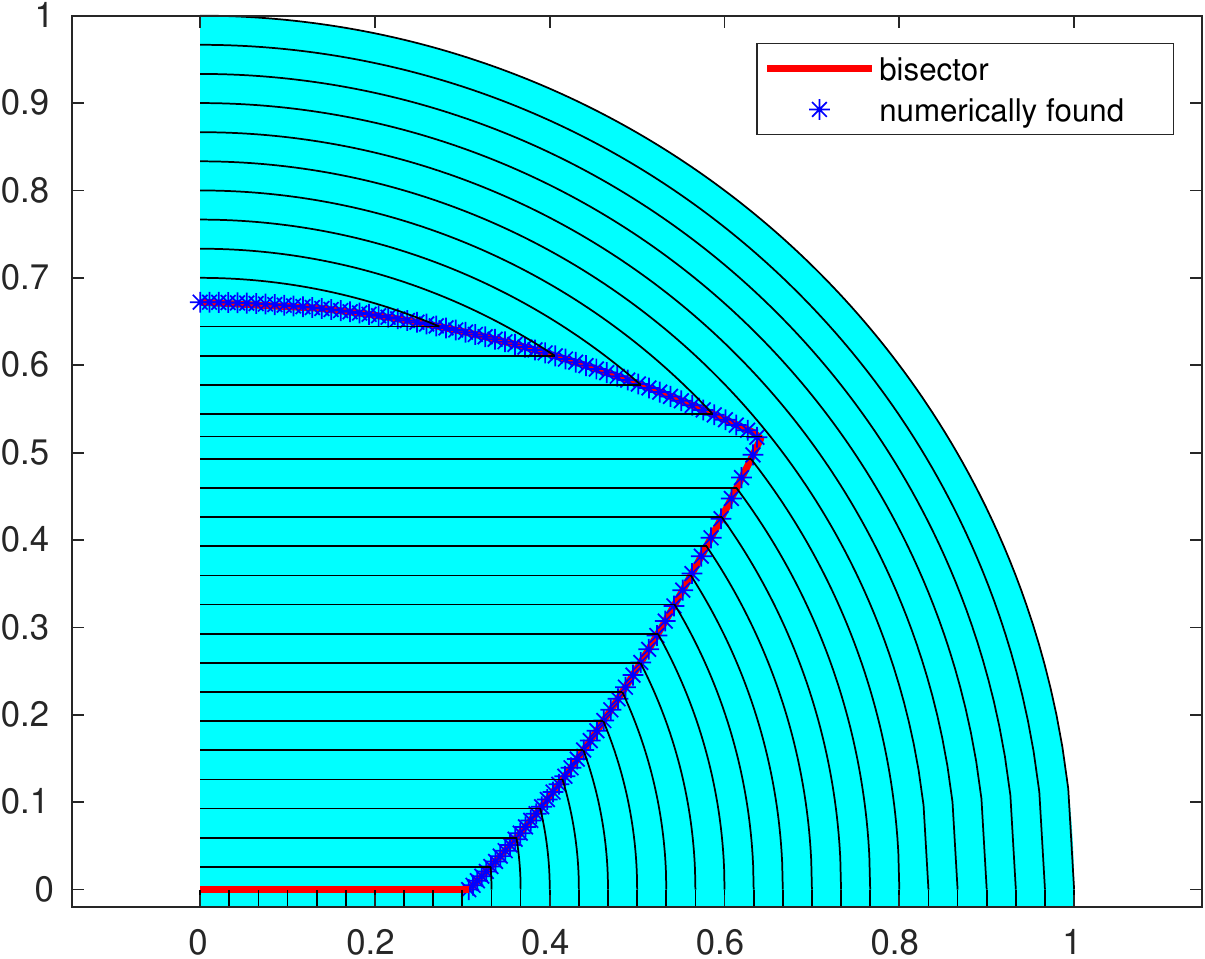}
   
   $\mu = 1$, $\alpha = 0.2$
  \end{minipage}
  
  \bigskip
  \begin{minipage}{.49\textwidth}
   \centering
   \includegraphics[width=.98\textwidth]{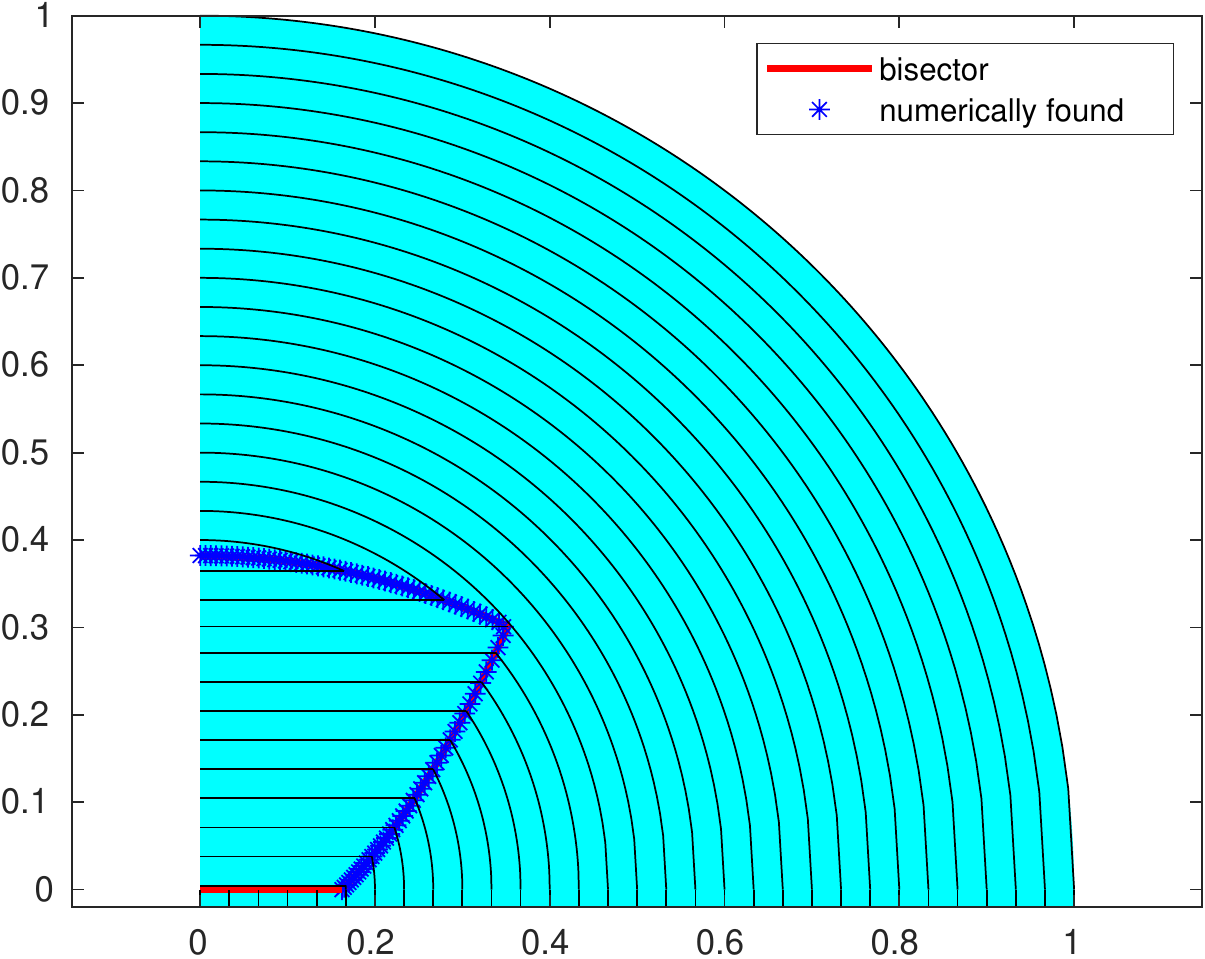}
   
   $\mu = 2$, $\alpha = 0.5$
  \end{minipage} \
  \begin{minipage}{.49\textwidth}
   \centering
   \includegraphics[width=.98\textwidth]{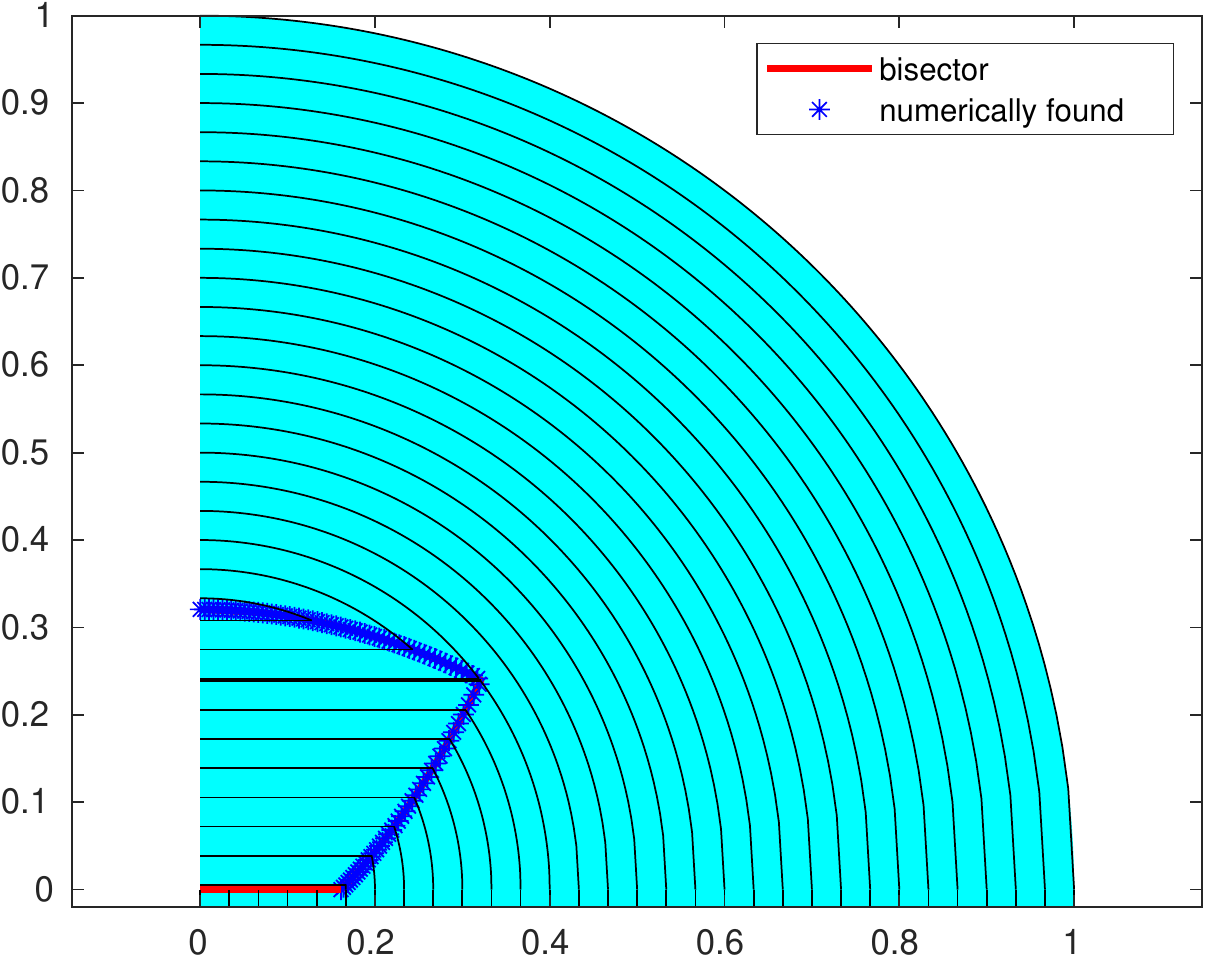}
   
   $\mu = 2$, $\alpha = 0.2$
  \end{minipage}
%
  \caption{Numerics for the problem on a quarter circle. All the pictures have~$K_1 = 2$.}
  \label{fig:numerics}
\end{figure}

Fig.~\ref{fig:numerics} shows a numerical approximation
of the minimizer of~\eqref{energy} within this restricted class, for~$K_1 = 2$
and different values of~$\alpha$, $\mu$. 
The simulation is based on a MATLAB code and the details
are provided in {{Appendix~\ref{app:numerics}}}.
The numerical method we use is \emph{not} guaranteed to
converge to a global minimum of the problem. However,
we repeated the simulations using several different initial guesses for~$\rho$,
including random ones, and obtained qualitatively
similar profiles for the jump set.
The numerically found jump set very nearly agrees with a union of
two parabolic arcs that bisect the smectic layers at each point,
given explicitly as
\begin{equation} \label{parabola-num}
 x_1 = \min\left(\left(a^2 + 2a x_2\right)^{1/2}, \,
 \left(b^2 - 2bx_2\right)^{1/2}\right)
\end{equation}
for some positive numbers~$a$, $b$. In Fig.~\ref{fig:numerics},
we compare the numerical solution with the
parabolic arcs given by~\eqref{parabola-num},
where we have chosen the parameters~$a$, $b$ so as to match the boundary values
of the numerically found solution; the difference between the two curves
is almost unnoticeable. Numerical tests show that the
jump set tends to shrink towards the centre of the circle
as~$\mu$ increases or~$\alpha$ decreses,
although the dependence on~$\alpha$ seems to be weaker.
This is consistent with the fact that the jump energy density
is monotonically increasing as a function of~$\mu$
and decreasing as a function of~$\alpha$.
Therefore, as~$\mu$ increases or~$\alpha$ decreases, 
energy minimization requires the jump set 
to become shorter, in order to compensate
for the larger values of the jump energy density.

\section{Discussion and future directions}

In this paper we have shown how a free-discontinuity model has the potential to describe the configurations of smectic A layers observed in thin film experiments. An advantage of free-discontinuity models in this context is that they represent defect walls as sharp interfaces with precise locations. However, our analysis is just a beginning and much remains to be done.  In particular, we are currently unable to formulate the thin film problem as a well-posed free-boundary problem in either two or three dimensions, so that the profile of the upper free surface and the width of the oily streaks can be predicted. 

Of course, it would be physically relevant 
to consider three-dimensional free-discontinuity problems.
The conditions~\eqref{C:first}--\eqref{C:rot}
are still meaningful in three dimensions, 
up to a few obvious modifications.
(As a side remark, assuming that the three-dimensional 
energy density is invariant with respect to the 
action of rotations~$\bar{\mathbf{R}}\in\mathrm{SO}(3)$ only
is enough to guarantee that, in the two-dimensional case,
the condition~\eqref{C:rot} is satisfied 
for any~$\mathbf{R}\in\mathrm{O}(2)$. Indeed,
any~$\mathbf{R}\in\mathrm{O}(2)$ can be realized as
a~$(2\times 2)$-submatrix of a three-dimensional
rotation~$\bar{\mathbf{R}}\in\mathrm{SO}(3)$.)
The definition of BV-ellipticity, too, extends to 
higher dimensions with no essential change 
(see e.g.~\cite[Definition~5.13]{AmbrosioFuscoPallara}
for the details). However, the condition~\eqref{C:disl}
needs more substantial modifications, because~\eqref{nodisl-n} 
(or equivalently, \eqref{nodisl}) 
is \emph{not} enough to avoid dislocations of the smectic layers at the
defect walls in three dimensions. In order for the layers
to be consistently defined across a defect wall,
the unit normal~$\nnu$ to the jump set 
and the molecular directors~$\n^+$, $\n^-$
on either side of the jump must satisfy~\eqref{nodisl-n}
and, additionally, they must be coplanar. 
(If~$\n^+$, $\n^-$ and~$\nnu$ are not coplanar,
then the intersection lines between the layers 
on either side do not belong to the jump set
and, hence, dislocations arise.)
The singular energy density~$\zeta_\alpha$, 
defined in~\eqref{zeta_alpha}, needs to be modified accordingly. 
Finally, the energy density~$\varphi_\alpha$ 
is still well-defined in three dimensions,
as~\eqref{phi_alpha} remains meaningful.
However, we do not know whether the three-dimensional
analogue of~$\varphi_\alpha$ is the BV-elliptic envelope
of the three-dimensional analogue of~$\zeta_\alpha$,
nor whether it is BV-elliptic at all, because 
our arguments in~\cite{BCS} do not apply 
to the three dimensional case.

Analytic work on these problems would need to be supplemented by numerical studies, and there is a need to develop appropriate numerical methods for free-discontinuity (and free boundary) problems associated with jump energy densities of the type~\eqref{phi_alpha}. Such jump energies also need further study, and it would be useful to broaden the class of  BV-elliptic jump energies satisfying \eqref{C:first}--\eqref{C:disl}, in particular to allow
 different minima at the two angle bisectors.
   
Finally, it is important to understand the relation between models based on a molecular density and our `sharp interface' model. In this context there are interesting numerical computations of 
 Xia, Maclachlan, Atherton \& Farrell~\cite{Xiaetal2021}
 using a modification of a model of Ball \& Bedford~\cite{BallBedford},
 which is in turn based on that of Pevnyi, Selinger \& Sluckin~\cite{PevnyiSelingerSluckin}.
 
\section*{Acknowledgement.} We would like to thank Emmanuelle Lacaze for many very helpful discussions. We also thank  Giacomo Albi and Marco Caliari from the University of Verona for their help with the numerical simulations, and the referee for their careful reading of the manuscript and comments.

\section*{Disclosure statement.} 
The authors report that there are no competing interests to declare.

\section*{Data availability statement.}
Data sharing is not applicable to this article as no new data were created or analyzed in this study.

 \newcommand{\noop}[1]{}

\section{Appendix: numerical simulations}
\label{app:numerics}

In this section, we present the details of the numerical
simulations for the minimization problem on a quarter circle,
described in Section~\ref{sect:results}. 
The $\Q$-tensor is defined in terms of a scalar 
function~$\rho\colon[0, \, \pi/2]\to\R$, as in~\eqref{Qrho}, 
\eqref{director-num}. However, we find it convenient to
introduce a new variable~$u = u(\theta)$, defined by
\begin{equation} \label{rho-u}
 \rho(\theta) = e^{-u(\theta)} 
 \qquad \textrm{for any } 0 < \theta < \frac{\pi}{2}.
\end{equation}
The original variable~$\rho$ is subject to the geometric 
constraints~$0 < \rho < 1$, because the domain is 
a quarter circle of unit radius.
However, we minimize numerically the energy functional
subject to \emph{no} constraints on~$u$. 
Equation~\eqref{rho-u} guarantees that~$\rho > 0$,
but for some values of~$\alpha$, $\mu$ (not shown in Fig.~\ref{fig:numerics}),
we did find numerical solutions that are not admissible, because
they do not satisfy~$\rho < 1$.

The energy can be written as a functional of~$u$
by a direct computation, using~\eqref{Qrho}, \eqref{director-num} 
and~\eqref{rho-u}. 
The energy consists of a sum of three terms:
\begin{equation} \label{Itot}
 I(u) = I^{\mathrm{el}}(u)
 + I^{\mathrm{jump,int}}(u)
 + I^{\mathrm{jump,bd}}(u)
\end{equation}
The first term, $I^{\mathrm{el}}(u)$, is the elastic energy:
\begin{equation} \label{Iel}
 I^{\mathrm{el}}(u) := \frac{K_1}{2}\int_{\Omega} \abs{\nabla\Q(\bfx)}^2 \, \d\bfx
 = \frac{K_1}{2} \int_0^{\pi/2} 
  \left( \int_{e^{-u(\theta)}}^1 \frac{\d\rho}{\rho} \right) \d\theta
 = \frac{K_1}{2} \int_0^{\pi/2} u(\theta) \, \d\theta
\end{equation}
The second term, $I^{\mathrm{jump,int}}(u)$,
is the contribution to the jump energy from 
the curve~$\mathcal{C}$, given in~\eqref{C-num}:
\begin{equation} \label{Ijump,int}
 \begin{split}
  I^{\mathrm{jump,int}}(u) 
  &:= \mu\int_{\mathcal{C}} \varphi(\Q^+, \, \Q^-, \, \nnu) \, \d\H^1 \\
  &= \mu\int_0^{\pi/2} \abs{\cos\theta}^\alpha \, e^{-u(\theta)}
   \left({u^\prime(\theta)}^2 + 1 + \abs{f(\theta, \, u^\prime(\theta))}\right)^{1/2} \d\theta
 \end{split}
\end{equation}
with
\begin{equation} \label{f}
 f(\theta, \, u^\prime(\theta)) := ({u^\prime(\theta)}^2 - 1)\cos\theta 
   + 2u^\prime(\theta)\,\sin\theta
\end{equation}
The expression~\eqref{Ijump,int}--\eqref{f} is deduced 
from~\eqref{phi_alpha_bis}, by applying trigonometric identities
to simplify the form of the integrand. The integrand  in~\eqref{Ijump,int}
is not differentiable at the points where~$f(\theta, \, u^\prime(\theta)) = 0$.
As the numerical minimization is based on a quasi-Newton method,
which works best for differentiable functions, we regularize the integrand by introducing a parameter~$\eps >0$:
\begin{equation} \label{Ijump,int-eps}
 \begin{split}
  I^{\mathrm{jump,int}}(u)
  \approx \mu\int_0^{\pi/2} \abs{\cos\theta}^\alpha \, e^{-u(\theta)}
   \left({u^\prime(\theta)}^2 + 1 + \left(\eps + f(\theta, \, u^\prime(\theta))^2\right)^{1/2} \right)^{1/2} \d\theta
 \end{split}
\end{equation}
The pictures in Fig.~\ref{fig:numerics} are obtained for~$\eps = 10^{-12}$.
Finally, $I^{\mathrm{jump,bd}}(u)$ is the contribution
to the jump energy from jumps that are located on
the boundary of~$\Omega$ --- more precisely, on the line 
segment~$S$ of endpoints~$(0, \, 0)$ and~$(0, \, \rho(0))$:
\begin{equation} \label{Ijump,bd}
 \begin{split}
  I^{\mathrm{jump,bd}}(u) 
  &:= \mu\int_{S} \varphi(\Q^+, \, \Q^-, \, \nnu) \, \d\H^1
  = \sqrt{2}\mu \, e^{-u(0)}
 \end{split}
\end{equation}
This term can also be written in integral form,
by considering a smooth function $g\colon[0,\,\pi/2]\to\R$
such that~$g(0) = 1$, $g(\pi/2) = 0$ and writing
\begin{equation} \label{Ijump,bd-bis}
 \begin{split}
  I^{\mathrm{jump,bd}}(u) 
  &= - \sqrt{2} \mu \int_0^{\pi/2} \frac{\d}{\d\theta}\left(e^{-u(\theta)} \, g(\theta)\right) \d\theta \\
  &= \sqrt{2} \mu \int_0^{\pi/2} e^{-u(\theta)}  \left(u^\prime(\theta) \,  g(\theta) - g^\prime(\theta)\right) \d\theta
 \end{split}
\end{equation}
We found that, when writing the boundary term 
in the form~\eqref{Ijump,bd}, the numerical solution
deviates from the parabolic arcs~\eqref{parabola-num}
near~$\theta = 0$. However, the thickness of this `boundary layer'
is mesh-dependent, so this feature is probably a numerical artifact.
The pictures in Fig.~\ref{fig:numerics} are obtained
by considering the integral form~\eqref{Ijump,bd-bis} of~$I^{\mathrm{jump,bd}}$,
with~$g(\theta) := 1 - 2\theta/\pi$, and present no boundary layer.
Other choices of the function~$g$, e.g.~$g(\theta) := \cos\theta$,
produce qualitatively similar profiles for the jump set.

We discretize the functional~\eqref{Itot} on
a uniform mesh of~$m$ points in~$[0, \, \pi/2]$.
We approximate the derivative~$u^\prime(\theta)$ 
by second-order central finite differences
and we approximate the integrals by the trapezoid rule.
The number~$m$ of mesh points is increased gradually,
from~$m = 50$ to~$m = 100$.
At each step of the iteration over~$m$,
we call the built-in MATLAB function \textsf{fminunc},
which applies the Broyden-Fletcher-Goldfarb-Shanno 
Quasi-Newton algorithm 
to minimize the functional~\eqref{Itot} under no constraints on~$u$.
The initial guess for the minimization process is defined
by the numerical minimizer found at the previous step.

\end{document}